\documentclass[12pt]{article}
\usepackage{amsmath,amssymb,amsthm,booktabs,fullpage}
\usepackage{latexsym,anyfontsize,mathtools,microtype}
\usepackage{aliascnt}
\usepackage[hidelinks]{hyperref}
\hypersetup{
    pdftitle={On Primitive Divisor Multiplicities of Lucas Sequences},
    pdfauthor={Samin Riasat},
    pdfsubject={Number theory}
}

\newtheorem{theorem}{Theorem}[section]

\newaliascnt{lemma}{theorem}
\newtheorem{lemma}[lemma]{Lemma}
\aliascntresetthe{lemma}

\newaliascnt{proposition}{theorem}
\newtheorem{proposition}[proposition]{Proposition}
\aliascntresetthe{proposition}

\newaliascnt{corollary}{theorem}
\newtheorem{corollary}[corollary]{Corollary}
\aliascntresetthe{corollary}

\theoremstyle{definition}
\newaliascnt{definition}{theorem}

\aliascntresetthe{definition}

\newaliascnt{question}{theorem}

\aliascntresetthe{question}

\theoremstyle{remark}
\newaliascnt{remark}{theorem}
\newtheorem{remark}[remark]{Remark}
\aliascntresetthe{remark}

\numberwithin{equation}{section}
\newcommand{\Z}{\mathbb{Z}}
\newcommand{\Q}{\mathbb{Q}}
\newcommand{\Fp}{\mathbb{F}}
\newcommand{\OK}{\mathcal{O}_{K}}

\newcommand{\rad}{\operatorname{rad}}
\newcommand{\vp}[1]{\nu_{#1}}
\newcommand{\eps}{\varepsilon}
\newcommand{\ph}{\varphi}
\newcommand{\mob}{\mu}
\newcommand{\tot}{\phi}
\newcommand{\leg}[2]{(#1 / #2)}
\newcommand{\Leg}[2]{\left(\frac{#1}{#2}\right)}
\allowdisplaybreaks[2]
\providecommand{\keywords}[1]{%
  \par\noindent
  {\small\textbf{\textit{Keywords---}} #1\par}%
}

\begin{document}
\title{On Primitive Divisor Multiplicities of Lucas Sequences}
\author{Samin Riasat}

\maketitle

\begin{abstract}
    Let $(U_{n})$ be a nondegenerate Lucas sequence with coprime integer parameters. Assuming the rational $abc$ conjecture, we prove that almost every index $n$ has a primitive prime divisor of $U_{n}$ of multiplicity one. More quantitatively, the logarithmic contribution of repeated primitive prime factors tends to zero relative to the primitive part along any family with $\phi(n) / n$ bounded below, and also in Ces\`aro mean over all indices. We derive explicit logarithmic lower bounds for Lucas non-Wieferich primes with ranks in prescribed residue classes. Both signs of the discriminant are allowed. Independently, we describe initial multiplicities at good odd primes by a unique Hensel root in the recurrence-parameter disk and by the exact contraction rate of a normalized Dickson-polynomial iteration. Applications include rank and period lifting for Fibonacci numbers and the Wall--Sun--Sun condition.
\end{abstract}


\section{Introduction}
Let $u, v \in \Z$ be coprime integers, $v \ne 0$, and let $\alpha, \beta$ be the roots of $X^{2} - uX + v$. We assume throughout that $\alpha / \beta$ is not a root of unity. In particular, $\Delta = u^{2} - 4v \ne 0$ and $u \ne 0$. Define the Lucas sequences of the first and second kinds by
\begin{align*}
    U_{n} = \frac{\alpha^{n} - \beta^{n}}{\alpha - \beta}, \qquad V_{n} = \alpha^{n} + \beta^{n} \qquad (n \ge 0).
\end{align*}
Equivalently, both satisfy 
\begin{align}
    \label{eq:rec}
    x_{n + 2} &= ux_{n + 1} - vx_{n}, 
\end{align}
with initial values $(x_{0}, x_{1}) = (0, 1)$ and $(2, u)$ respectively. Put $A = \max(|\alpha|, |\beta|) > 1$. For a prime $p \nmid v$, let $k_{U}(p)$ denote its rank of apparition, the least positive $n$ with $p \mid U_{n}$, and write $r_{U}(p) = \vp{p}(U_{k_{U}(p)})$, where $\vp{p}$ is the $p$-adic valuation, with $\vp{p}(0) = +\infty$. Such a rank always exists (see \autoref{prop:rank}). We call $p$ a primitive divisor at index $n$ if $k_{U}(p) = n$, and simple if additionally $r_{U}(p) = 1$. The usual Lucas primitive-divisor convention excludes primes dividing $\Delta$; our first-occurrence convention includes them. This affects only finitely many indices for any fixed sequence.

The classical primitive-divisor theorem ensures that, for $n > 30$, $U_{n}$ has a primitive divisor in the usual sense \cite{BHV}; for real roots, Carmichael's earlier theorem gives a stronger threshold \cite{Carmichael}. For a comprehensive account of Lucas sequences, see Ballot and Williams \cite{BW}. These existence results do not control the multiplicity of a primitive divisor. For the Fibonacci sequence $(u, v) = (1, -1)$, the assertion that every primitive divisor is simple is equivalent to the nonexistence of Wall--Sun--Sun primes. In general it is false: even within $\Q(\sqrt{5})$, one can construct Lucas sequences in which a prescribed odd prime $p \ne 5$ has initial multiplicity at least two (\autoref{prop:subseq}).

We investigate two aspects of initial multiplicity: its global distribution over the indices, conditional on $abc$, and its local description as a function of the recurrence parameter. For the first, write
\begin{equation}
    \label{eq:SH}
    P_{n} = \prod_{k_{U}(p) = n}p^{r_{U}(p)}, \qquad
    S_{n} = \prod_{\substack{k_{U}(p) = n\\
    r_{U}(p) = 1}}p, \qquad
    H_{n} = \prod_{\substack{k_{U}(p) = n\\
    r_{U}(p) \ge 2}}p^{r_{U}(p)}.
\end{equation}
Thus $P_{n} = S_{n}H_{n}$. Empty products are $1$. Throughout, $\tot$ denotes Euler's totient function, $\mob$ the M\"obius function, and $\omega(n)$ the number of distinct prime factors of $n$. By the rational $abc$ conjecture we mean that, for each $\epsilon > 0$, coprime positive integers $a + b = c$ satisfy $c \ll_{\epsilon}\rad(abc)^{1 + \epsilon}$, where $\rad(m) = \prod_{p \mid m}p$ for nonzero $m$.

\begin{theorem}[Conditional abundance of simple primitive divisors]
    \label{thm:abundance}
    Assume the rational $abc$ conjecture. For every fixed nondegenerate Lucas sequence with coprime integer parameters the following hold.
    \begin{enumerate}
        \item For every $\eta, \delta > 0$, uniformly over indices $n > 1$ with $\tot(n) \ge \eta n$,
        \begin{equation}
            \label{eq:Sbound}
            H_{n} \ll_{u, v, \eta, \delta}P_{n}^{\delta},
            \qquad S_{n} \gg_{u, v, \eta, \delta}P_{n}^{1 - \delta}.
        \end{equation}
        In particular, $\log S_{n} / \log P_{n} \to 1$ along any unbounded family of such indices.
        \item Put $R_{n} = \log H_{n} / \log P_{n}$ if $P_{n} > 1$, and $R_{n} = 0$ otherwise. Then
        \begin{align*}
            \lim_{X \to \infty}\frac{1}{X}\sum_{n \le X}R_{n}
            &= 0.
        \end{align*}
        \item The indices $n$ with $S_{n} > 1$ have natural density one. Moreover, for every fixed prime $p$, all sufficiently large indices $p^{a}$ have $S_{p^{a}} > 1$.
    \end{enumerate}
    These conclusions also hold if primitive divisors are required to be coprime to the discriminant.
\end{theorem}

This result refines the conditional method of extracting non-Wieferich primes from cyclotomic factors. Ribenboim and Walsh \cite{RW} obtained bounds on powerful parts of binary recurrent sequences of positive discriminant; see also Yabuta \cite{Yabuta} for negative discriminant and Ribenboim \cite{Ribenboim} for conditional infinitude of Lucas non-Wieferich primes. Peng \cite{Peng} treats a generalization of Fibonacci--Wieferich primes. Most directly related to our argument is the work of Anitha, Fathima and Vijayalakshmi \cite{AFV}, who, using the number-field $abc$ conjecture, obtain logarithmic lower bounds for Lucas non-Wieferich primes in the residue classes $\pm 1$ modulo a fixed integer. Their cyclotomic extraction uses a positive-proportion family of indices with $\tot(n) / n$ bounded below. We use the rational Lucas identity to work with rational $abc$ and obtain a density-one statement for ranks themselves, without a restriction on the sign of the discriminant. The powerful-part bound and cyclotomic estimates are established tools; the focus of \autoref{thm:abundance} is their consequence for the distribution of simple primitive divisors.

\begin{corollary}[Conditional rank distribution]
    \label{cor:rankdistribution}
    Assume rational $abc$. For $h \ge 1$ and $b \pmod{h}$, let $N_{h, b}(x)$ count primes $p \le x$, $p \nmid 2v\Delta$, with $r_{U}(p) = 1$ and $k_{U}(p) \equiv b \pmod{h}$. Then
    \begin{equation}
        \label{eq:count}
        \liminf_{x \to \infty}\frac{N_{h, b}(x)}{\log x} \ge \frac{1}{h\log A}.
    \end{equation}
    In particular, for $h \ge 2$ the same lower bound holds for the count of Lucas non-Wieferich primes $p \le x$ with $p \equiv \pm 1 \pmod{h}$.
\end{corollary}
Here a prime $p \nmid 2v\Delta$ is Lucas non-Wieferich if $\vp{p}(U_{p - \leg{\Delta}{p}}) = 1$, equivalently $r_{U}(p) = 1$ by the rank and lifting formulas, 
where $\leg{\cdot}{p}$ is the Legendre symbol.
The conclusion for arbitrary $b$ restricts the rank rather than just the prime's residue class.

Our second result concerns variation of $u$ at fixed $v$. Write $U_{n}(X, v), V_{n}(X, v) \in \Z[X]$ for the Lucas polynomials obtained by replacing $u$ by $X$ in \eqref{eq:rec}. They are Dickson polynomials of the second and first kinds, respectively (with index $n - 1$ in the second-kind case); see \cite{Dickson}. At an odd prime $p \nmid v\Delta$, set
\begin{equation}
    \label{eq:Tdef}
    \lambda_{p} = \Leg{v}{p}\, v^{-(p - 1) / 2} \in \Z_{p}^{\times},
    \qquad T_{p}(X) = \lambda_{p} V_{p}(X, v).
\end{equation}
The normalization ensures $\lambda_{p}^{2}v^{p} = v$ and $\lambda_{p} \equiv 1 \pmod{p}$, so that the iteration preserves the norm parameter and the residue class of $u$.

\begin{theorem}[Local multiplicity and normalized iteration]
    \label{thm:hensel}
    Let $p \nmid 2v\Delta$ be prime, $k = k_{U}(p)$ and $r = r_{U}(p)$. There is a unique $\widehat{u} \in u + p\Z_{p}$ with $U_{k}(\widehat{u}, v) = 0$. For $x, y \in u + p\Z_{p}$,
    \begin{align}
        \vp{p}(U_{k}(x, v))
        &= \vp{p}(x - \widehat{u}),
        \label{eq:distance}\\
        \vp{p}(T_{p}(x) - T_{p}(y))
        &= 1 + \vp{p}(x - y).
        \label{eq:contract}
    \end{align}
    In particular, $\widehat{u}$ is the unique fixed point of $T_{p}$ in this disk. Setting $u_{j} = T_{p}^{j}(u)$, for every $j \ge 0$ one has
    \begin{align}
        u_{j}
        &= \lambda_{p}^{(p^{j} - 1) / (p - 1)}V_{p^{j}}(u, v),
        \label{eq:iter}\\
        \vp{p}(u_{j} - \widehat{u})
        &= \vp{p}(u_{j + 1} - u_{j}) = r + j.
        \label{eq:iterdist}
    \end{align}
    Consequently $\vp{p}(T_{p}(u) - u) = r$, and the exceptional parameter modulo $p^{a}$ is uniquely determined for every $a \ge 1$.
\end{theorem}

For $v = -1$, the normalization is trivial, and the iteration consists of the usual trace polynomials. In particular, for Fibonacci numbers $T_{p}^{j}(1) = L_{p^{j}}$ and $\vp{p}(L_{p} - 1) = r_{U}(p)$. Bragman and Rowland \cite[Theorem 10]{BR} express the Fibonacci Wall exponent using the principal-unit part of the golden ratio. The local description here is consistent with theirs: its focus is the recurrence parameter, the norm-preserving normalization for arbitrary $v$, and the exact valuation of each iteration step. McIntosh and Roettger \cite[Theorem 1(v)]{McR} already give a congruence characterization of Lucas--Wieferich primes in terms of $V_{p}(u, v)$ for $p \nmid 2uv\Delta$. Our valuation identity refines this characterization to all exponents and also allows $p \mid u$. The underlying Dickson composition and power-map dynamics are classical \cite{Dickson,DS}. For related congruences and derivative identities for Lucas polynomials, see Yamagishi \cite{Yamagishi}; M\"uller \cite{Muller} studies the distribution of ranks when a recurrence parameter varies over a finite field.

The proofs are organized as follows. \autoref{sec:lift} records the lifting and rank results needed below, including the discriminant-prime case. \autoref{sec:cyclo} establishes the primitive-factor estimates, including the prime $2$. \autoref{sec:abc} proves \autoref{thm:abundance} and \autoref{cor:rankdistribution}; \autoref{sec:padic} proves \autoref{thm:hensel}. Finally, \autoref{sec:fib} gives Fibonacci applications, drawing on the author's earlier notes \cite{Rperiod,Rprime,Rblog,Rlte}.

\section{Lifting in rings of algebraic integers}
\label{sec:lift}
Throughout, $\vp{p}$ is the usual $p$-adic valuation on integers, with $\vp{p}(0) = +\infty$. We also use $\vp{p}$ for the valuation on $\Q_{p}$ normalized by $\vp{p}(p) = 1$. For a number field $K$ we use congruence modulo $p^{j}$ in $\OK$ to mean congruence modulo the ideal $p^{j}\OK$. In particular, an integer in $p^{j}\OK$ belongs to $p^{j}\Z$, since $\OK\cap\Q = \Z$.

\begin{theorem}[Algebraic lifting]
    \label{thm:alg}
    Let $p$ be an odd rational prime, and let $\alpha, \beta \in \OK$ be units modulo $p\OK$. Put
    \begin{align*}
        A_{n}
        &= \begin{cases}
            \sum_{i = 0}^{n - 1}\alpha^{n - 1 - i}\beta^{i}, & n \ge 1, \\ 
            0, & n = 0. 
        \end{cases}
    \end{align*}
    Suppose $0 \ne A_{k} \in p\OK$. For every positive integer $m$ and prime ideal $\mathfrak{p} \mid p$,
    \begin{equation}
        \label{eq:ideal}
        v_{\mathfrak{p}}(A_{km}) = v_{\mathfrak{p}}(A_{k}) + e_{\mathfrak{p}}\vp{p}(m),
    \end{equation}
    where $e_{\mathfrak{p}} = v_{\mathfrak{p}}(p)$ and $v_{\mathfrak{p}}$ is the normalized ideal valuation. In particular, if $A_{n}$ is integral over $\Z$ and rational for every $n$, then
    \begin{equation}
        \label{eq:lift}
        \vp{p}(A_{km}) = \vp{p}(A_{k}) + \vp{p}(m).
    \end{equation}
\end{theorem}
\begin{proof}
    The polynomial identity $A_{kt} = A_{k} S_{t}(c, d)$ holds with $c = \alpha^{k}$, $d = \beta^{k}$ and $S_{t}(c, d) = \sum_{i = 0}^{t - 1}c^{t - 1 - i}d^{i}$. Since $c - d = (\alpha - \beta)A_{k} \in p\OK$, we have $c \equiv d \pmod{p}$. If $p \nmid t$, then $S_{t}(c, d) \equiv t d^{t - 1} \pmod{p}$ is a unit modulo $p\OK$.

    For $t = p$, write $c = d + pz$, where $z \in \OK$. The binomial theorem gives
    \begin{align*}
        S_{p}(c, d)
        &= \sum_{j = 1}^{p}\binom{p}{j} d^{p - j}(pz)^{j - 1}
        \equiv pd^{p - 1} \pmod{p^{2}}.
    \end{align*}
    Here the intermediate terms are in $p^{2}\OK$, and so is the last term because $p \ge 3$. Hence $S_{p}(c, d) = p w$ with $w$ a unit modulo $p\OK$. Replacing $c, d$ by their $p$th powers preserves their congruence and unit properties. Factor $m = p^{a} t$ with $p \nmid t$ and apply the two calculations successively. The factors prime to every $\mathfrak{p} \mid p$ have valuation zero, and each factor $p$ has valuation $e_{\mathfrak{p}}$. This proves \eqref{eq:ideal}; for rational integers it implies \eqref{eq:lift}.
\end{proof}
\begin{remark}
    No unramifiedness hypothesis is required. For $p = 2$ the same step is valid if $c \equiv d \pmod{4}$, since $S_{2}(c, d) = c + d \equiv 2d \pmod{4}$. It need not hold when only $c \equiv d \pmod{2}$; for instance, $F_{3} = 2$ but $F_{6} = 8$.
\end{remark}

\subsection{Ranks of apparition}
\begin{proposition}
    \label{prop:rank}
    The sequence $(U_{n})$ satisfies $\gcd(U_{m}, U_{n}) = |U_{\gcd(m, n)}|$ for positive $m, n$. For every $M > 1$ coprime to $v$, the least positive index $k_{U}(M)$ with $M \mid U_{k_{U}(M)}$ exists, and
    \begin{align*}
        M \mid U_{n}\quad
        &\Longleftrightarrow \quad k_{U}(M) \mid n.
    \end{align*}
    For an odd prime $p \nmid v$, put $\eps = \leg{\Delta}{p}$. Then
    \begin{equation}
        \label{eq:rank-bound}
        k_{U}(p) \mid p - \eps, \qquad U_{p} \equiv \eps \pmod{p}.
    \end{equation}
    If $p \mid \Delta$, then $k_{U}(p) = p$. A prime dividing $v$ divides no $U_{n}$ with $n \ge 1$.
\end{proposition}
\begin{proof}
    The identity $U_{m + n} = U_{m}U_{n + 1} - vU_{m - 1}U_{n}$ follows from the recurrence \eqref{eq:rec}. Every prime dividing $v$ is coprime to each positive-index $U_{n}$, since $U_{n} \equiv u^{n - 1}$ modulo that prime. The recurrence then shows $\gcd(U_{m}, U_{m - 1}) = 1$ by descending induction. The addition identity implies $\gcd(U_{m}, U_{m + n}) = \gcd(U_{m}, U_{n})$, so the Euclidean algorithm proves strong divisibility. The recurrence state $(U_{n}, U_{n + 1})$ evolves by an invertible matrix modulo $M$ if $\gcd(M, v) = 1$. It returns to $(0, 1)$, giving existence; strong divisibility and minimality yield the equivalence.

    Modulo $p$, if $\Delta$ is a nonzero square, the two nonzero roots lie in $\Fp_{p}$ and have $(p - 1)$st powers equal to $1$. If $\Delta$ is a nonsquare, the roots lie in $\Fp_{p^{2}}$ and Frobenius interchanges them, so their $(p + 1)$st powers are both $v$. In either case the Binet expression gives $U_{p - \eps} = 0$ and $U_{p} = \eps$ modulo $p$. In the remaining case the polynomial has a double root $c = u / 2 \ne 0$ in $\Fp_{p}$; the polynomial sum defining $U_{n}$ reduces to $nc^{n - 1}$. This gives $k_{U}(p) = p$ and $U_{p} \equiv 0 \pmod{p}$.
\end{proof}

\begin{corollary}
    \label{cor:ranklift}
    Let $p$ be an odd prime not dividing $v$, $k = k_{U}(p)$ and $r = \vp{p}(U_{k})$. Then for every positive $n$ and $a$,
    \begin{align}
        \vp{p}(U_{n})
        &= \begin{cases}
            r + \vp{p}(n / k), &k \mid n, \\
            0, &k \nmid n,
        \end{cases}
        \label{eq:valuation}\\
        k_{U}(p^{a})
        &= p^{\max(0, a - r)}k.
        \label{eq:ranklift}
    \end{align}
\end{corollary}
\begin{proof}
    Apply \autoref{thm:alg} to the roots of $X^{2} - uX + v$ (both are units modulo $p$ because their product is $v$) and use \autoref{prop:rank}. The smallest multiplier with valuation at least $a - r$ gives \eqref{eq:ranklift}.
\end{proof}

\subsection{Primes dividing the discriminant}
\begin{proposition}
    \label{prop:disc}
    For odd primes $p$ one has $p \mid U_{p}$ if and only if $p \mid \Delta$. If $p \mid \Delta$ and $p > 3$, then
    \begin{equation}
        \label{eq:disc}
        U_{p} \equiv p \pmod{p^{2}}.
    \end{equation}
    Moreover $9 \mid U_{3}$ if and only if $\Delta \equiv 6 \pmod{9}$.
\end{proposition}
\begin{proof}
    For $p \nmid v$ the first assertion follows from \autoref{prop:rank}; if $p \mid v$, then $p \nmid U_{p}\Delta$ by coprimality. The symmetric expansion is
    \begin{align*}
        2^{p - 1}U_{p}
        &= \sum_{j = 0}^{(p - 1) / 2}\binom{p}{2j + 1}u^{p - 2j - 1}\Delta^{j}.
    \end{align*}
    If $p \mid \Delta$ and $p > 3$, every term except $pu^{p - 1}$ is divisible by $p^{2}$, including the last term since $(p - 1) / 2 \ge 2$. Since $p \nmid u$, Fermat's little theorem gives \eqref{eq:disc}. Finally $4U_{3} = 3u^{2} + \Delta$. If $9 \mid U_{3}$ then $3 \mid \Delta$, and $u^{2} \equiv 1 \pmod{3}$. Conversely $\Delta \equiv 6 \pmod{9}$ implies $3 \nmid u$; the displayed identity proves the equivalence.
\end{proof}
\begin{remark}
    This corrects the one-sided binomial expansion in \cite[Lemma 4.1]{Rlte}: $p \mid \Delta$ does not in general imply $\sqrt{\Delta} \in p\OK$. The symmetric expansion avoids that issue. The valuation formulas agree with the general results of Sanna \cite{Sanna}. In particular, $\vp{p}(U_{n}) = \vp{p}(n)$ for $p > 3$ dividing $\Delta$; for $p = 3 \mid \Delta$ and $3 \mid n$ it is $\vp{3}(U_{3}) + \vp{3}(n) - 1$.
\end{remark}
\section{Cyclotomic factors and their primitive parts}
\label{sec:cyclo}
For $n > 1$ let
\begin{equation}
    \label{eq:cyclo}
    C_{n} = |\Phi_{n}(\alpha, \beta)| = \prod_{d \mid n}|U_{d}|^{\mob(n / d)},
\end{equation}
where $\Phi_{n}(X, Y) = Y^{\tot(n)}\Phi_{n}(X / Y)$ is the homogeneous cyclotomic polynomial. The identity follows by M\"obius inversion from the cyclotomic factorization of $\alpha^{n} - \beta^{n}$. The polynomial is symmetric for $n > 1$, so its value is a rational algebraic integer and $C_{n}$ is a positive integer. Moreover $|U_{n}| = \prod_{d \mid n, \, d > 1}C_{d}$. These factors, also known as Lucas atoms or M\"obius duals, have been studied in \cite{ST,AMMR,RSC}. We record their valuations to control the nonprimitive contribution, with a short proof from \autoref{thm:alg}.

\begin{proposition}[{\cite{AMMR,RSC}}]
    \label{thm:cyclo}
    Let $n > 1$, let $p \nmid v$ be prime, and write $k = k_{U}(p)$, $r = r_{U}(p)$. If $p$ is odd, then
    \begin{equation}
        \label{eq:Cval}
        \vp{p}(C_{n}) = \begin{cases}
            r, &n = k, \\
            1, &n = p^{a}k, \ a \ge 1, \\
            0, &\text{otherwise}.
        \end{cases}
    \end{equation}
    If $p = 2$, then $v$ is odd, $k = 2$ for even $u$ and $k = 3$ for odd $u$, and
    \begin{equation}
        \label{eq:C2}
        \vp{2}(C_{n}) = \begin{cases}
            r, &n = k, \\
            \vp{2}(V_{k}), &n = 2k, \\
            1, &n = 2^{a}k, \ a \ge 2, \\
            0, &\text{otherwise}.
        \end{cases}
    \end{equation}
    A prime dividing $v$ divides no $C_{n}$.
\end{proposition}
\begin{proof}
    For odd $p$, insert \eqref{eq:valuation} into the M\"obius transform in \eqref{eq:cyclo}. If $k \nmid n$, all terms vanish. If $n = km$, the transform is
    \begin{equation*}
        \sum_{t \mid m}\mob(m / t)(r + \vp{p}(t)).
    \end{equation*}
    The constant contributes only at $m = 1$. The transform of $\vp{p}(t) = \sum_{a \ge 1}\mathbf{1}_{p^{a} \mid t}$ is $1$ precisely at $m = p^{a}$ for $a \ge 1$. This also covers primes dividing $\Delta$, with their own initial value $r$.

    For $p = 2$, the ranks follow from the recurrence \eqref{eq:rec} modulo $2$. If $2 \mid U_{j}$, the roots $\alpha^{j}, \beta^{j}$ are congruent modulo $2$ in their ring of integers and are units modulo $2$, so the quotient $U_{jt} / U_{j}$ is odd for odd $t$ by the geometric-sum identity. Since $V_{j} = 2U_{j + 1} - uU_{j}$, $V_{k}$ is even. The identities
    \begin{align*}
        U_{2j}
        &= U_{j} V_{j}, \qquad V_{2j} = V_{j}^{2} - 2v^{j}
    \end{align*}
    show that the first doubling increases the valuation by $\vp{2}(V_{k})$, and each subsequent doubling increases it by one, since $V_{2j} \equiv 2 \pmod{4}$ whenever $V_{j}$ is even and $v$ is odd. Thus $\vp{2}(U_{kt}) = r$ for odd $t$, and $r + \vp{2}(V_{k}) + \vp{2}(t) - 1$ for even $t$. Its M\"obius transform gives \eqref{eq:C2}. Finally, $p \mid v$ implies $p \nmid U_{d}$ for every positive $d$ by \autoref{prop:rank}, hence $p \nmid C_{n}$.
\end{proof}

\begin{corollary}
    \label{cor:D}
    For $n > 1$, write $C_{n} = D_{n}P_{n}$. The odd part of $D_{n}$ is either $1$ or the largest prime factor $Q$ of $n$, occurring to exponent one. In particular,
    \begin{equation}
        \label{eq:Dgrowth}
        1 \le D_{n} \ll_{u, v}n.
    \end{equation}
\end{corollary}
\begin{proof}
    By \autoref{thm:cyclo}, a nonprimitive odd divisor $p$ requires $n = p^{a}k_{U}(p)$ for $a \ge 1$. If $p \mid \Delta$, then $k_{U}(p) = p$; otherwise $k_{U}(p) \mid p - \leg{\Delta}{p}$, whose prime factors are all smaller than $p$ (also for $p + 1$, which is even). Thus $p = Q$ and its exponent is one. The contribution of $2$ is bounded by $2^{\max(1, \vp{2}(V_{k_{U}(2)}))}$ if $v$ is odd, and is absent otherwise. This proves the bound.
\end{proof}

\begin{lemma}[Growth of the primitive part]
    \label{lem:growth}
    For every fixed nondegenerate Lucas sequence with coprime integer parameters,
    \begin{equation}
        \label{eq:growth}
        \log P_{n} \sim \tot(n)\log A\qquad(n \to \infty).
    \end{equation}
    If $\Delta > 0$, the sharper estimate $\log P_{n} = \tot(n)\log A + O_{u, v}(\log n)$ holds.
\end{lemma}
\begin{proof}
    Choose $|\alpha| = A \ge |\beta|$. Since $\alpha\beta = v$ is a nonzero integer, $A \ge 1$. Equality would force both absolute values to be one; if the roots are real their quotient is then a root of unity, and if complex, $v = 1$ and $u \in \{-1, 0, 1\}$, again giving a root of unity. Hence $A > 1$. By \eqref{eq:cyclo},
    \begin{equation}
        \label{eq:logcyclo}
        \log C_{n} = \tot(n)\log A +
        \sum_{d \mid n}\mob(n / d)\log|1 - (\beta / \alpha)^{d}|.
    \end{equation}
    For real roots, $|\beta / \alpha| < 1$ by nondegeneracy. The absolute sum of the error terms over all $d \ge 1$ converges, so $\log C_{n} = \tot(n)\log A + O_{u, v}(1)$.

    For complex conjugate roots the ratio $z = \beta / \alpha$ has modulus one but is not a root of unity. The standard bound for a linear form in logarithms, in the form of \cite[Lemma 3]{Voutier}, gives that the logarithm $\log|1 - z^{d}|$ is at least $-c_{u, v}(\log(2d))^{2}$ for some positive constant; it is at most $\log2$. Therefore the error in \eqref{eq:logcyclo} is $O_{u, v}(2^{\omega(n)}(\log(2n))^{2})$. The elementary estimates $2^{\omega(n)} = n^{o(1)}$ and $\tot(n) = n^{1 - o(1)}$ show this is $o(\tot(n))$. For completeness, for each $\epsilon > 0$ the inequalities $2 \le p^{\epsilon}$ and $1 - 1 / p \ge p^{-\epsilon}$ hold for all sufficiently large primes $p$; the remaining primes contribute a fixed constant to the two products. Combining with \eqref{eq:Dgrowth} proves both assertions.
\end{proof}

\section{Conditional density of simple primitive divisors}
\label{sec:abc}
We use the rational Lucas identity to bound the entire powerful part, and then remove the nonprimitive factors. The estimate below is a version of the powerful-part bounds of Ribenboim and Walsh \cite{RW}; the proof also works for negative discriminant, as in the setting of Yabuta \cite{Yabuta}.

\begin{lemma}[Conditional powerful-part bound]
    \label{thm:abc}
    Assume rational $abc$. For $n \ge 1$, define the powerful part
    \begin{align*}
        T_{n}
        &= \prod_{\vp{p}(U_{n}) \ge 2}p^{\vp{p}(U_{n})}.
    \end{align*}
    For each $\delta > 0$, it satisfies $T_{n} \ll_{u, v, \delta}A^{\delta n}$. In particular, $\log H_{n} = o(n)$.
\end{lemma}
\begin{proof}
    The identity
    \begin{equation}
        \label{eq:LucasPell}
        V_{n}^{2} - \Delta U_{n}^{2} = 4v^{n}
    \end{equation}
    has three nonzero terms: nondegeneracy excludes $U_{n} = 0$ and $V_{n} = 0$ for $n \ge 1$. If $p \mid v$, then $U_{n} \equiv u^{n - 1}$ and $V_{n} \equiv u^{n} \pmod{p}$, so neither is divisible by $p$. It follows that the common divisor $g$ of the three terms in \eqref{eq:LucasPell} divides $4$ (if $v$ is even then $u$ is odd and the common divisor is $1$). Dividing by $g$ gives pairwise coprime integers, and moving the negative terms to the other side gives an equation of the form $a + b = c$ with $a, b, c > 0$.

    The maximum of the absolute values of the original terms is $\asymp_{u, v} A^{2n}$. For real roots, this follows from $|V_{n}| \asymp_{u, v} A^{n}$; for complex roots, use $|4v^{n}| = 4A^{2n}$. Also $|U_{n}|, |V_{n}| \ll_{u, v} A^{n}$. The radical of $4v^{n}\Delta$ is fixed, so the conjecture implies, for every $\epsilon > 0$,
    \begin{align*}
        A^{2n}
        &\ll_{u, v, \epsilon}\bigl(A^{n}\rad(|U_{n}|)\bigr)^{1 + \epsilon}.
    \end{align*}
    The radical is thus $\gg_{u, v, \epsilon} A^{n(1 - \epsilon) / (1 + \epsilon)}$. On the other hand,
    \begin{align*}
        \rad(|U_{n}|)
        &= \frac{|U_{n}|}{T_{n}}\rad(T_{n}) \le \frac{|U_{n}|}{\sqrt{T_{n}}}
        \ll_{u, v}\frac{A^{n}}{\sqrt{T_{n}}}.
    \end{align*}
    Combining gives $T_{n} \ll_{u, v, \epsilon} A^{4\epsilon n / (1 + \epsilon)}$. Taking $\epsilon$ sufficiently small proves the claim. Since $H_{n} \mid T_{n}$, the final assertion follows.
\end{proof}

\begin{proof}[Proof of \autoref{thm:abundance}]
    By \autoref{lem:growth}, for $n$ sufficiently large we have $\log P_{n} \ge \tfrac{1}{2}\tot(n)\log A$. In particular, if $\tot(n) \ge \eta n$, apply \autoref{thm:abc} with exponent $\eta\delta / 2$ to obtain $H_{n} \ll_{u, v, \eta, \delta}P_{n}^{\delta}$. The bound for $S_{n} = P_{n} / H_{n}$ follows. The finitely many small indices can be absorbed in the constants. Taking arbitrary $\delta > 0$ proves the logarithmic limiting assertion.

    For the mean-value assertion, observe that
    \begin{equation}
        \label{eq:mean}
        \frac{1}{X}\sum_{n \le X}\frac{n}{\tot(n)} \le
        \sum_{d \ge 1}\frac{\mob(d)^{2}}{d\tot(d)}
        = \prod_{p\text{ prime}}\left(1 + \frac{1}{p(p - 1)}\right) \eqqcolon C < \infty.
    \end{equation}
    Indeed, the identity $n / \tot(n) = \sum_{d \mid n}\mob(d)^{2} / \tot(d)$ follows by multiplicativity, and the product converges by comparison with $\sum_{m \ge 2}m^{-2}$. By \autoref{thm:abc}, for each $\epsilon > 0$, $\log H_{n} \le \epsilon n$ for all sufficiently large $n$. Hence for such $n$,
    \begin{align*}
        R_{n}
        &\le \frac{2\epsilon}{\log A}\frac{n}{\tot(n)}.
    \end{align*}
    The finite initial portion contributes $o(1)$ to the mean. By \eqref{eq:mean} the limsup of the mean is at most $2C\epsilon / \log A$, for every $\epsilon > 0$, and is therefore zero.

    By \autoref{lem:growth}, $P_{n} > 1$ for all sufficiently large $n$. When $S_{n} = 1$ and $P_{n} > 1$ we have $R_{n} = 1$. The mean-value conclusion therefore shows that the set of such indices has density zero. For $n = p^{a}$ we have $\tot(n) / n = 1 - 1 / p$, so the first assertion with any $\delta < 1$ and \autoref{lem:growth} show that $S_{n} > 1$ for all sufficiently large $a$. Finally, removing the primes dividing $\Delta$ alters only the finitely many indices at which such a prime first appears; the estimates and limits remain valid after a change of constants.
\end{proof}

\begin{proof}[Proof of \autoref{cor:rankdistribution}]
    By \autoref{thm:abundance}, the number of indices $n \le Y$ in the class $b \pmod{h}$ with $S_{n} > 1$ is $Y / h + o(Y)$. Excluding the finitely many ranks of primes dividing $2v\Delta$, choose a simple primitive divisor $p_{n}$ of each such $U_{n}$. The primes are distinct because the ranks are distinct. Since $p_{n} \le |U_{n}| \le c_{u, v}A^{n}$, taking $Y = \lfloor\log(x / c_{u, v}) / \log A\rfloor$ proves \eqref{eq:count}. Finally, when $b = 0$ the rank bound in \autoref{prop:rank} gives $h \mid p_{n} - \leg{\Delta}{p_{n}}$, so $p_{n} \equiv \pm 1 \pmod{h}$. The equivalence of simplicity with the stated Lucas non-Wieferich condition follows from \eqref{eq:valuation} because $p_{n}$ is coprime to the index $p_{n} - \leg{\Delta}{p_{n}}$.
\end{proof}

\begin{remark}
    These statements concern the distribution of ranks and not a prime-density asymptotic. They neither assert that every sufficiently large index has a simple primitive divisor nor exclude repeated primitive divisors. The argument uses the unspecified constants of the rational $abc$ conjecture and gives no explicit rate for its exceptional set.
\end{remark}

\section{Local multiplicity in recurrence-parameter space}
\label{sec:padic}
We prove \autoref{thm:hensel}. In this section $v$ is held fixed. The Binet expressions for the Lucas polynomials give the derivative identities (the second is also \cite[Lemma 2.1]{Yamagishi})
\begin{equation}
    \label{eq:deriv}
    (X^{2} - 4v)\frac{\partial U_{n}}{\partial X} = nV_{n} - XU_{n},
    \qquad \frac{\partial V_{n}}{\partial X} = nU_{n}.
\end{equation}
Indeed, if $a, b$ are the roots of $Y^{2} - XY + v$ over $\Q(X)$, then $a' = a / (a - b)$ and $b' = -b / (a - b)$. Differentiating the Binet expressions gives \eqref{eq:deriv}. We also have $V_{p}(X, v) \equiv X^{p} \pmod{p}$ as polynomials by the Frobenius identity for the trace; see also \cite[Lemma 3.1]{Yamagishi}.

\begin{proof}[Proof of \autoref{thm:hensel}]
    The rank bound yields $k \mid p - \leg{\Delta}{p}$ and $p \nmid k$. By \eqref{eq:LucasPell}, $p \nmid V_{k}(u, v)$ because $p \mid U_{k}(u, v)$ and $p \nmid 2v$. The first identity in \eqref{eq:deriv} gives $U_{k}'(u, v) \not\equiv 0 \pmod{p}$. Hensel's lemma supplies a unique $\widehat{u}$ in $u + p\Z_{p}$. On that disk the polynomial divided difference of $U_{k}(X, v)$ at any two points is congruent to $U_{k}'(u, v)$ modulo $p$, hence a unit. Taking one point to be $\widehat{u}$ proves \eqref{eq:distance}.

    Euler's criterion implies $\lambda_{p} \equiv 1 \pmod{p}$; the normalization also gives $\lambda_{p}^{2}v^{p} = v$. Consequently $T_{p}(x) \equiv x \pmod{p}$ for all $x \in \Z_{p}$. For $x$ in the disk, the derivative $T_{p}'(x) = \lambda_{p} p U_{p}(x, v)$ has valuation one, by \eqref{eq:rank-bound} applied at the residue of $u$. To obtain the exact contraction, write $y = x + h$, $s = \vp{p}(h) \ge 1$, and expand the polynomial $T_{p}(x + h)$. Its coefficient of $h$ has valuation one. The coefficients of $h^{i}$ for $2 \le i < p$ are in $p\Z_{p}$, since $T_{p}(X) \equiv X^{p} \pmod{p}$; the coefficient of $h^{p}$ is a unit. All these higher terms have valuation strictly greater than $1 + s$ (here $p \ge 3$). The linear term therefore gives \eqref{eq:contract}. If $x = y$, both sides have the usual infinite-valuation interpretation.

    Let $c_{j} = \lambda_{p}^{(p^{j} - 1) / (p - 1)}$. The roots at parameter $u_{j} = T_{p}^{j}(u)$ are $c_{j}\alpha^{p^{j}}$ and $c_{j}\beta^{p^{j}}$. This follows by induction, as their product is $c_{j}^{2} v^{p^{j}} = v$ and the next iteration scales their $p$th powers by $\lambda_{p}$. It proves \eqref{eq:iter}. The Binet expression gives
    \begin{align*}
        U_{k}(u_{j}, v)
        &= c_{j}^{k - 1}\frac{U_{kp^{j}}(u, v)}{U_{p^{j}}(u, v)}.
    \end{align*}
    Here $c_{j}$ is a $p$-adic unit and the denominator is a unit, since $k > 1$ and $p \nmid k$. By \eqref{eq:valuation}, the valuation of the numerator is $r + j$. Applying \eqref{eq:distance} yields $\vp{p}(u_{j} - \widehat{u}) = r + j$. Thus $u_{j} \to \widehat{u}$, and continuity gives $T_{p}(\widehat{u}) = \widehat{u}$. The two successive distances to $\widehat{u}$ have unequal valuations, so $\vp{p}(u_{j + 1} - u_{j}) = r + j$. Finally, \eqref{eq:contract} rules out two distinct fixed points. This proves the theorem.
\end{proof}

\begin{corollary}
    \label{cor:exception}
    Under the hypotheses of \autoref{thm:hensel}, for each $a \ge 1$ there is exactly one residue class modulo $p^{a}$ in $u + p\Z_{p}$ on which $p^{a} \mid U_{k}(x, v)$. Of its $p$ lifts modulo $p^{a + 1}$, exactly $p - 1$ have valuation $a$ and one has higher valuation. In particular,
    \begin{align*}
        \vp{p}(U_{k}(u, v))
        &= \vp{p}(\lambda_{p} V_{p}(u, v) - u), \qquad
        \widehat{u} \equiv T_{p}(u) \pmod{p^{r + 1}}.
    \end{align*}
\end{corollary}
\begin{proof}
    The class is $x \equiv \widehat{u} \pmod{p^{a}}$ by \eqref{eq:distance}; the remaining statements are immediate from \eqref{eq:iterdist}.
\end{proof}
To compare the valuation identity with \cite[Theorem 1(v)]{McR}, put $w = \leg{v}{p}\, v^{(p - 1) / 2} = \lambda_{p}^{-1}$. Since $w \equiv 1 \pmod{p}$, we have $w^{2} + 1 \equiv 2w \pmod{p^{2}}$. Thus $\lambda_{p}^{-1} \equiv (v^{p - 1} + 1) / 2 \pmod{p^{2}}$, and $\lambda_{p} V_{p}(u, v) \equiv u \pmod{p^{2}}$ is precisely their trace congruence. The exceptional class at precision $p^{2}$ has relative frequency $1 / p$ as the parameter varies in the disk. This exact statement about a variable parameter does not imply equidistribution over primes for a fixed integer $u$.

For example, take $(u, v) = (3, 2)$ and $p = 7$. Here $U_{n} = 2^{n} - 1$, $V_{n} = 2^{n} + 1$, $k = 3$, $r = 1$, and $\lambda_{7} = 1 / 8$. Hence $T_{7}(3) = 129 / 8$ and $\vp{7}(T_{7}(3) - 3) = \vp{7}(105 / 8) = 1$, as predicted. The first iterate satisfies $T_{7}(3) \equiv 10 \pmod{49}$, the unique lift of $3$ modulo $7$ solving $U_{3}(X, 2) = X^{2} - 2 \equiv 0 \pmod{49}$.

The Dickson normalization in \eqref{eq:Tdef} allows any fixed nonzero norm parameter, rather than just $v = -1$. The iteration itself is a rescaled power map, with classical Dickson composition; our application records the rank exponent as its exact initial distance and step size. For $v = -1$, $\lambda_{p} = 1$, so no rescaling is needed.

\section{Fibonacci applications}
\label{sec:fib}
\subsection{Ranks and periods}
Specialize to $(u, v) = (1, -1)$, writing $F_{n} = U_{n}$, $L_{n} = V_{n}$, $\ph = (1 + \sqrt{5}) / 2$, and $\psi = (1 - \sqrt{5}) / 2 = -\ph^{-1}$. Write $k(m) = k_{U}(m)$ for the rank of apparition modulo $m$ and $\ell(m)$ for the Pisano period modulo $m$. Put $r_{p} = \vp{p}(F_{k(p)})$ for each prime $p$, and $\eps_{p} = \leg{5}{p}$ for odd $p \ne 5$. The quantities $C_{n}, D_{n}, P_{n}, S_{n}, H_{n}$ in this section are their Fibonacci specializations. A prime $p \ne 2, 5$ is Wall--Sun--Sun when $\vp{p}(F_{p - \eps_{p}}) \ge 2$. Since $r_{2} = r_{5} = 1$ and $r_{p} = \vp{p}(F_{p - \eps_{p}})$ otherwise, the assertion that every first-occurring prime is simple is equivalent to the nonexistence of Wall--Sun--Sun primes. The problem goes back to Wall \cite{Wall}; see also Sun and Sun \cite{SunSun}. No such prime is known; the customary heuristic instead predicts infinitely many \cite{McR,BR}.

In $R = \Z[\ph]$ the identity
\begin{equation}
    \label{eq:phi}
    \ph^{n} = F_{n}\ph + F_{n - 1}\qquad(n \ge 1)
\end{equation}
follows by induction. The ring $R$ is the ring of integers of $\Q(\sqrt{5})$. The period $\ell(M)$ is exactly the order of $\ph$ in $(R / MR)^{\times}$, by \eqref{eq:phi} and the independence of $1, \ph$ as a $\Z$-basis. This interpretation uses ideal congruence and does not depend on whether a prime dividing $M$ splits, ramifies or remains inert in $\Q(\sqrt{5})$.

\begin{proposition}
    \label{prop:period}
    For $M > 1$, $\ell(M) / k(M)$ divides $4$. If $p$ is an odd prime, then
    \begin{equation}
        \label{eq:periodlift}
        \ell(p^{a}) = p^{\max(0, a - r_{p})}\ell(p)\qquad(a \ge 1).
    \end{equation}
\end{proposition}
\begin{proof}
    Let $k = k(M)$. Equation \eqref{eq:phi} implies $\ph^{k} \equiv F_{k - 1} \pmod{M}$. Taking norms gives $F_{k - 1}^{2} \equiv (-1)^{k} \pmod{M}$. Thus $F_{k - 1}$ is a unit of order dividing $4$. The powers of $\ph$ that are rational scalars modulo $M$ occur exactly at multiples of $k$, by strong divisibility. Hence $\ell(M) / k(M)$ equals the order of $F_{k - 1}$ modulo $M$ and divides $4$.

    The kernel of reduction $(R / p^{a}R)^{\times} \to (R / pR)^{\times}$ has order a power of $p$ (it is $1 + pR / p^{a}R$). Therefore $\ell(p^{a}) / \ell(p)$ is a power of $p$. On the other hand $k(p^{a}) / k(p)$ is a power of $p$ by \eqref{eq:ranklift}, and the ratios of period to rank are powers of $2$. Since $p$ is odd, these latter ratios are equal at $p^{a}$ and $p$. This gives \eqref{eq:periodlift}.
\end{proof}

\begin{corollary}
    \label{cor:wss}
    For an odd prime $p \ne 5$, the following are equivalent: $p$ is Wall--Sun--Sun; $r_{p} \ge 2$; $k(p^{2}) = k(p)$; $\ell(p^{2}) = \ell(p)$.
\end{corollary}
\begin{proof}
    Use $k(p) \mid p - \eps_{p}$ from \eqref{eq:rank-bound}, so $p \nmid (p - \eps_{p}) / k(p)$, and apply \eqref{eq:valuation}, \eqref{eq:ranklift} and \eqref{eq:periodlift}.
\end{proof}

\begin{theorem}[Odd-rank congruences]
    \label{thm:odd-rank}
    Let $p \ne 2, 5$ be prime with odd rank $n = k(p)$. Then $\ell(p) = 4n$ and
    \begin{equation}
        \label{eq:oddcong}
        p \equiv \begin{cases}
            1 \pmod{4n}, &\eps_{p} = 1, \\
            2n - 1 \pmod{4n}, &\eps_{p} = -1.
        \end{cases}
    \end{equation}
    In particular, $p \equiv 1 \pmod{4}$.
\end{theorem}
\begin{proof}
    As in \autoref{prop:period}, $F_{n - 1}^{2} \equiv -1 \pmod{p}$; for odd $p$ its order is $4$. Thus $\ell(p) = 4n$. In the split case $\ph^{p - 1} \equiv 1 \pmod{p}$ in $R / pR \cong \Fp_{p} \times \Fp_{p}$, giving $4n \mid p - 1$. In the inert case Frobenius sends $\ph$ to $\psi = -\ph^{-1}$, so $\ph^{p + 1} \equiv -1 \pmod{p}$. Since $\ph^{2n} \equiv -1 \pmod{p}$, we have $4n \mid p + 1 - 2n$. This proves the stated congruences.
\end{proof}

\begin{corollary}
    \label{cor:prime-index}
    If $q > 5$ is prime and $p$ is a prime divisor of $F_{q}$, then $k(p) = q$, $\ell(p) = 4q$, and \eqref{eq:oddcong} holds with $n = q$. If $e = \vp{p}(F_{q})$, then $\ell(p^{a}) = 4p^{\max(0, a - e)}q$ for every $a \ge 1$. In particular, $p^{2} \mid F_{q}$ implies that $p$ is a Wall--Sun--Sun prime.
\end{corollary}
\begin{proof}
    The rank divides $q$ and exceeds $1$. Also $p \ne 2, 5$ since their ranks are $3, 5$. Apply \autoref{thm:odd-rank} and \autoref{cor:wss}.
\end{proof}
\begin{remark}
    This extends the prime-index classification of \cite[Proposition 2]{Rprime} to all odd ranks. It does not imply that a hypothetical Wall--Sun--Sun prime has prime, or even odd, rank. In particular, squarefreeness of every prime-index Fibonacci number is only a necessary consequence of nonexistence of Wall--Sun--Sun primes.
\end{remark}

\subsection{Cyclotomic factors and splitting}
For Fibonacci numbers, $k(2) = 3$, $r_{2} = 1$, and $\vp{2}(L_{3}) = 2$. Thus \eqref{eq:C2} has the exceptional value $2$ at $n = 6$ and is at most one elsewhere. The odd part of $D_{n}$ is either $1$ or the largest prime factor of $n$, as in \autoref{cor:D}. We recall the squarefreeness characterization of Ross, Shen and Cai.
\begin{corollary}[{\cite[Corollary 5.1]{RSC}}]
    \label{cor:squarefree}
    The following assertions are equivalent:
    \begin{enumerate}
        \item No Wall--Sun--Sun prime exists.
        \item $P_{n}$ is squarefree for every $n > 1$.
        \item The odd part of $C_{n}$ is squarefree for every $n > 1$.
    \end{enumerate}
\end{corollary}
\begin{proof}
    The first two are equivalent by definition and $r_{2} = r_{5} = 1$. For an odd prime $p$, \eqref{eq:Cval} shows that $p^{2} \mid C_{n}$ if and only if $k(p) = n$ and $r_{p} \ge 2$. This proves the third equivalence. By \eqref{eq:C2}, the third condition is also equivalent to squarefreeness of $C_{n}$ at every $n \ne 6$. This is exactly the characterization of Ross, Shen and Cai \cite[Corollary 5.1]{RSC}.
\end{proof}

\subsubsection*{Splitting and parity of primitive multiplicities}
The valuation formula above is a specialization of \cite[Theorems 14--15]{AMMR} and \cite[Theorem 2.2]{RSC}. Ross, Shen and Cai \cite[Corollaries 2.5--2.6]{RSCcong} establish the parity of primitive multiplicities for general regular Lucas sequences and explicitly for Fibonacci numbers. In particular, their Corollary 2.6(b), (e) includes the inert prime-power existence result below. We provide a direct proof using the character modulo $5$ and combine it with the odd-rank congruences to obtain the modulus $4n$. Their Corollary 2.4(b) also includes the prime-power quotient congruence modulo $p^{a}$.
We next restrict to $n > 1$ with $\gcd(n, 30) = 1$, so every prime dividing $C_{n}$ is odd and different from $5$. This restriction is convenient rather than essential. For odd $d$ one has
\begin{equation}
    \label{eq:mod5}
    F_{d} \equiv d\, 3^{d - 1} \pmod{5},
    \qquad \Leg{F_{d}}{5} = \Leg{d}{5}.
\end{equation}
Indeed, the polynomial $T^{2} - T - 1$ has a double root $3$ modulo $5$. Its divided difference $F_{d}$ reduces to $d3^{d - 1}$, and $3^{d - 1}$ is a square when $d$ is odd. (The first congruence holds also for even $d$.)

\begin{theorem}[Character of the primitive part]
    \label{thm:char}
    Let $n > 1$ and $\gcd(n, 30) = 1$. Define
    \begin{align*}
        \eta_{n}
        &= \begin{cases}
            \leg{p}{5}, &n = p^{a}\text{ for a prime }p, \ a \ge 1, \\
            1, &n\text{ is not a prime power}.
        \end{cases}
    \end{align*}
    Then
    \begin{align}
        \Leg{C_{n}}{5}
        &= \eta_{n},
        \label{eq:Cchar}\\
        \Leg{P_{n}}{5}
        &= \eta_{n}\Leg{D_{n}}{5},
        \label{eq:Pchar}\\
        P_{n}
        &\equiv \begin{cases}
            1 \pmod{4n}, &\leg{P_{n}}{5} = 1, \\
            2n - 1 \pmod{4n}, &\leg{P_{n}}{5} = -1.
        \end{cases}
        \label{eq:Pcong}
    \end{align}
    In particular, the parity of the total multiplicity of inert primitive divisors of $F_{n}$ is determined by $n$ and $D_{n}$.
\end{theorem}
\begin{proof}
    By \eqref{eq:cyclo} and \eqref{eq:mod5},
    \begin{align*}
        \Leg{C_{n}}{5}
        &= \prod_{d \mid n}\Leg{d}{5}^{\mob(n / d)}
        = \Leg{\prod_{d \mid n}d^{\mob(n / d)}}{5}.
    \end{align*}
    The integer inside the symbol in the final expression equals $p$ for $n = p^{a}$ and $1$ otherwise, by taking the valuation at each prime and using the same M\"obius transform as in the proof of \autoref{thm:cyclo}. This proves \eqref{eq:Cchar}, and the identity $C_{n} = D_{n}P_{n}$ proves \eqref{eq:Pchar}.

    If $k(p) = n$, \autoref{thm:odd-rank} implies that $p \equiv 1$ or $2n - 1$ modulo $4n$ according as $p$ splits or is inert. The latter residue is an involution since $(2n - 1)^{2} \equiv 1 \pmod{4n}$. Consequently $P_{n}$ is congruent to the latter residue precisely when the sum of the multiplicities of inert factors is odd. By quadratic reciprocity, $\leg{5}{p} = \leg{p}{5}$, so that parity is measured by $\leg{P_{n}}{5}$. This proves \eqref{eq:Pcong}.
\end{proof}

\begin{corollary}[Prime-power indices]
    \label{cor:pp}
    Let $p > 5$ be prime and $a \ge 1$. Put $B_{p, a} = F_{p^{a}} / F_{p^{a - 1}}$. Then
    \begin{equation}
        \label{eq:B}
        B_{p, a} = C_{p^{a}} = P_{p^{a}}
        \equiv \begin{cases}
            1 \pmod{4p^{a}}, &\leg{p}{5} = 1, \\
            2p^{a} - 1 \pmod{4p^{a}}, &\leg{p}{5} = -1.
        \end{cases}
    \end{equation}
    If $p \equiv \pm 2 \pmod{5}$, then $F_{p^{a}}$ has an inert primitive divisor $q$ of odd multiplicity. Such a $q$ satisfies $q \equiv 2p^{a} - 1 \pmod{4p^{a}}$ and $q \equiv 13$ or $17 \pmod{20}$.
\end{corollary}
\begin{proof}
    For $n = p^{a}$, \eqref{eq:cyclo} reduces to the quotient in the statement. If $D_{n} \ne 1$, \autoref{cor:D} would give $n = p^{b}k(p)$ for $b \ge 1$. But $1 < k(p) \mid p - \eps_{p}$ and $p \nmid k(p)$, a contradiction. Thus $D_{n} = 1$ and \eqref{eq:B} follows from \autoref{thm:char}. When $\leg{p}{5} = -1$, the total inert multiplicity is odd and at least one such prime has odd multiplicity. The congruences follow from \autoref{thm:odd-rank} and quadratic reciprocity.
\end{proof}

\begin{remark}
    For a fixed inert prime $p > 5$, \autoref{cor:pp} gives infinitely many distinct inert primes, one of odd primitive multiplicity at each rank $p^{a}$. Each such multiplicity is either one or at least three; the character argument cannot distinguish the two possibilities. For example, with $p = 7$, $a = 2$, one has
    \begin{align*}
        B_{7, 2}
        &= 598364773 = 97 \cdot 6168709.
    \end{align*}
    The prime $97$ is inert, has rank $49$, and occurs exactly once. This is an illustration, not a general proof of simple primitive divisibility.
\end{remark}

\subsection{The local Wall exponent}
For Fibonacci numbers, $T_{p}(X) = V_{p}(X, -1)$ by \eqref{eq:Tdef}, so \autoref{thm:hensel} gives an exact description in terms of Lucas numbers.
\begin{corollary}
    \label{cor:Lucas}
    For every odd prime $p \ne 5$,
    \begin{equation}
        \label{eq:Lq}
        \vp{p}(L_{p} - 1) = r_{p} = \vp{p}(F_{p - \eps_{p}}).
    \end{equation}
    Moreover, $L_{p^{j}} \to \widehat{u}$ in $\Z_{p}$ as $j \to \infty$, and $\vp{p}(L_{p^{j}} - \widehat{u}) = r_{p} + j$ for every $j \ge 0$. In particular, $p$ is Wall--Sun--Sun if and only if $L_{p} \equiv 1 \pmod{p^{2}}$.
\end{corollary}
\begin{proof}
    Set $(u, v) = (1, -1)$ in \autoref{thm:hensel} and apply \eqref{eq:valuation}.
\end{proof}
The congruence characterization at exponent two is classical \cite[Theorem 1(v)]{McR}. The fixed point is also explicit. In the unramified local extension containing $\ph$, let $\zeta$ be its Teichm\"uller representative. Then the representative of $\psi = -\ph^{-1}$ is $-\zeta^{-1}$ and
\begin{equation}
    \label{eq:teich}
    \widehat{u} = \zeta - \zeta^{-1}.
\end{equation}
The sum belongs to $\Z_{p}$ (in the inert case Frobenius interchanges the two summands); their quotient has order $k(p)$, and their difference is a unit as $p \ne 5$. Thus $U_{k(p)}(\zeta - \zeta^{-1}, -1) = 0$, and Hensel uniqueness proves \eqref{eq:teich}. In fact
\begin{align*}
    1 - \widehat{u}
    &= (\ph - \zeta)(1 + (\ph\zeta)^{-1}),
\end{align*}
where the second factor is a unit, as it reduces to $\sqrt{5} / \ph$ modulo $p$. This recovers the local Wall exponent of Bragman and Rowland \cite[Theorem 10]{BR}.

\begin{proposition}
    \label{prop:subseq}
    Let $k \ge 1$. If $p > 3$ is a prime and $p \mid F_{k}$, then
    \begin{align*}
        \frac{F_{kp}}{F_{k}}
        &\equiv p \pmod{p^{2}}.
    \end{align*}
    Independently, for every odd prime $p \ne 5$, the Lucas sequence $B_{m} = F_{pm} / F_{p}$ has discriminant coprime to $p$, norm parameter $v = -1$, and first occurrence of $p$ at $m = k(p)$ with multiplicity $r_{p} + 1$.
\end{proposition}
\begin{proof}
    The first quotient is the $p$th term of the sequence $A_{m} = F_{km} / F_{k}$, whose parameters are $u = L_{k}$, $v = (-1)^{k}$, and $\Delta = 5F_{k}^{2}$. Thus \autoref{prop:disc} applies. For the second assertion use the same construction with index $p$. Now $p \nmid F_{p}$ by \autoref{prop:rank}, and $p \nmid k(p)$. Hence $p \mid B_{m}$ if and only if $k(p) \mid m$. Its initial valuation is $\vp{p}(F_{pk(p)}) = r_{p} + 1$ by \eqref{eq:valuation}.
\end{proof}
This example shows that a general multiplicity-one theorem for Lucas sequences is false, even within the field $\Q(\sqrt{5})$ and with norm parameter $-1$. Any proof for Fibonacci numbers must use additional information about their fixed recurrence parameter.

\section{Further questions}
The conditional density result allows repeated primitive divisors. For Fibonacci numbers it does not imply that $H_{n} = 1$ at every index, which would exclude all Wall--Sun--Sun primes. An unconditional result establishing $S_{p^{a}} > 1$ for infinitely many $a$ for a single fixed prime $p$ would already prove the infinitude of non-Wall--Sun--Sun primes with prescribed prime-power ranks. A related question, suggested by \autoref{cor:pp}, is whether for some inert prime $p > 5$ infinitely many $F_{p^{a}}$ have an inert simple primitive divisor. The quadratic-character constraint alone only guarantees an inert divisor of odd multiplicity, not of multiplicity one.

Locally, the multiplicity-one conjecture for Fibonacci numbers is the assertion $\vp{p}(1 - \widehat{u}) = 1$ for every $p \ne 2, 5$. The exact contraction of \autoref{thm:hensel} controls all further lifts but not that initial distance. The subsequence construction of \autoref{prop:subseq} illustrates the obstruction to a parameter-independent multiplicity-one assertion. Neither the conditional global estimate nor the local reformulation resolves this open problem.

\section*{Declaration of generative AI use}

The author used ChatGPT (OpenAI, GPT-6)
to assist in the consolidation and presentation of
previously developed mathematical results, the
exploration of further mathematical consequences,
literature searches, proof refinement, and manuscript
preparation. The research builds on the author's
earlier independent work on Fibonacci and Lucas
sequences. 
The author takes full responsibility
for the originality, accuracy, and integrity of
the contents.

\bibliographystyle{plain}
\bibliography{references}

\medskip
\noindent MSC2020: 11B39, 11A07, 11A41, 11D75, 11S05

\end{document}